\documentclass[11pt,a4paper]{amsart}
\usepackage{graphicx,multirow,array,amsmath,amssymb}

\newtheorem{theorem}{Theorem}
\newtheorem{lemma}[theorem]{Lemma}

\begin{document}

\title{Exactly fourteen intrinsically knotted graphs have 21 edges}
\author[M. Lee]{Minjung Lee} 
\address{Department of Mathematics, Korea University, Anam-dong, Sungbuk-ku, Seoul 136-701, Korea}
\email{mjmlj@korea.ac.kr}
\author[H. Kim]{Hyoungjun Kim}
\address{Department of Mathematics, Korea University, Anam-dong, Sungbuk-ku, Seoul 136-701, Korea}
\email{kimhjun@korea.ac.kr}
\author[H. J. Lee]{Hwa Jeong Lee}
\address{Department of Mathematical Sciences, KAIST, 291 Daehak-ro, Yuseong-gu, Daejeon 305-701, Korea}
\email{hjwith@kaist.ac.kr}
\author[S. Oh]{Seungsang Oh}
\address{Department of Mathematics, Korea University, Anam-dong, Sungbuk-ku, Seoul 136-701, Korea}
\email{seungsang@korea.ac.kr}

\thanks{The corresponding author(Seungsang Oh) was supported by the National Research Foundation of 
Korea(NRF) grant funded by the Korea government(MSIP) (No. NRF-2014R1A2A1A11050999).}
\thanks{The third author was supported by the National Research Foundation of Korea Grant
funded by the Korean Government (NRF-2010-0024630).}

\begin{abstract}
Johnson, Kidwell, and Michael showed that intrinsically knotted graphs have at least 21 edges.
Also it is known that $K_7$ and the thirteen graphs obtained from $K_7$ by $\nabla Y$ moves are
intrinsically knotted graphs with 21 edges.
We prove that these 14 graphs are the only intrinsically knotted graphs with 21 edges.
\end{abstract}

\maketitle

\section{Introduction} \label{sec:intro}

Throughout the article we will take an embedded graph to mean a graph embedded in $R^3$.
We call a graph $G$ {\em intrinsically knotted\/} if every embedding of the graph contains a knotted cycle.
Conway and Gordon~\cite{CG} showed that $K_7$, the complete graph with seven vertices,
is an intrinsically knotted graph.
A graph $H$ is {\em minor\/} of another graph $G$ if it can be obtained from $G$ by contracting or deleting some edges.
An intrinsically knotted graph is {\em minor minimal intrinsically knotted\/}
provided no proper minor is intrinsically knotted.
Robertson and Seymour~\cite{RS} proved that there are only finite minor minimal intrinsically knotted graphs,
but finding the complete set of them is still an open problem.
However, it is well known that $K_7$ and the thirteen graphs obtained from this graph by $\nabla Y$ moves are
minor minimal intrinsically knotted; see Conway and Gordon \cite{CG}, and Kohara and Suzuki \cite{KS}.

\begin{figure}[h]
\includegraphics{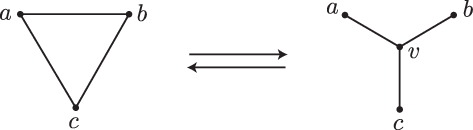}
\caption{$\nabla Y$ and $Y \nabla$ moves}
\label{fig:1}
\end{figure}

A $\nabla Y$ {\em move\/} is an exchanging operation that removes all edges of a triangle $abc$ and
inserts a new vertex $v$ and three edges $va, vb$ and $vc$ as in Figure~\ref{fig:1}.
Its reverse operation is called a $Y \nabla$ {\em move\/}.
Since $\nabla Y$ move preserves intrinsical knottedness 
(see Motwani, Raghunathan, and Saran \cite{MRS}),
we will only consider triangle-free graphs in the article.

From the work of Johnson, Kidwell, and Michael~\cite{JKM}
it follows that any intrinsically knotted graph consists at least 21 edges.
Here is the main theorem.

\begin{theorem}\label{thm:main1}
The only triangle-free intrinsically knotted graphs with exactly 21 edges are $H_{12}$ and $C_{14}$.
$($$H_{12}$ and $C_{14}$ were described by Kohara and Suzuki in \cite{KS}.$)$
\end{theorem}

\begin{figure}[h]
\includegraphics{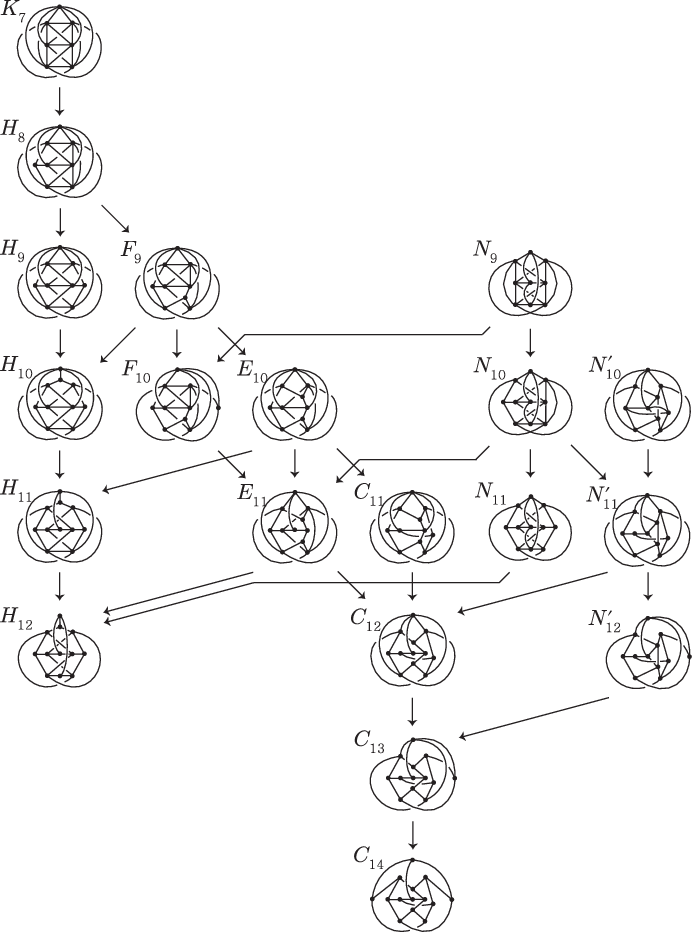}
\caption{The graph $K_7$ and $19$ more related graphs, where each arrow represents a $\nabla Y$ move}
\label{fig:2}
\end{figure}

Kohara and Suzuki \cite{KS} found fourteen intrinsically knotted graphs.
Goldberg, Mattman, and Naimi~\cite{GMN} constructed twenty graphs
derived from $H_{12}$ and $C_{14}$ by $Y \nabla$ moves as in Figure~\ref{fig:2},
and they showed that these six graphs $N_9$, $N_{10}$, $N_{11}$, $N'_{10}$, $N'_{11}$, and
$N'_{12}$ are not intrinsically knotted.
This fact was proved by Hanaki, Nikkuni, Taniyama, and Yamazaki~\cite{HNTY} independently.
Theorem \ref{thm:main1} guarantees that all intrinsically knotted graphs with 21 edges
can be obtained from $H_{12}$ and $C_{14}$ by $Y \nabla$ moves.
Thus we have the following theorem.

\begin{theorem}\label{thm:main2}
The only intrinsically knotted graphs with exactly 21 edges are
$K_7$ and the thirteen graphs obtained from $K_7$ by $\nabla Y$ moves.
\end{theorem}

This theorem gives us the complete set of fourteen minor minimal intrinsically knotted graphs with 21 edges.

\section{Terminology} \label{sec:term}

From now on let $G=(V,E)$ denote a triangle-free graph with $21$ edges.
Here $V$ and $E$ denote the sets of all vertices and edges of $G$, respectively.
For any two distinct vertices $a$ and $b$,
let $\widehat{G}_{a,b}=(\widehat{V}_{a,b}, \widehat{E}_{a,b})$ denote the graph
obtained from $G$ by deleting two vertices $a$ and $b$,
and then contracting an edge incident to a vertex of
degree $1$ or $2$ repeatedly until no vertices of degree $1$ or $2$ exist.
Removing vertices means deleting interiors of all edges incident to these vertices
as well as the resulting isolated vertices.

In a graph, the distance between two vertices $a$ and $b$ is the number of edges
in the shortest path connecting them and is denoted by dist$(a,b)$.
The degree of a vertex $a$ is denoted by deg$(a)$.
To count the number of edges of $\widehat{G}_{a,b}$, we introduce some notation.

\begin{itemize}
\item $E(a)$ is the set of edges which are incident to $a$
\item $V(a)=\{c \in V\ |\ \mbox{dist}(a,c)=1\}$
\item $V_n(a)=\{c \in V\ |\ \mbox{dist}(a,c)=1,\ \mbox{deg}(c)=n,\}$
\item $V_n(a,b)=V_n(a) \cap V_n(b)$
\item $V_Y(a,b)=\{c \in V\ |\ \exists \ d \in V_3(a,b) \ \mbox{such that} \ c \in V_3(d) \setminus \{a,b\}\}$
\end{itemize}

\begin{figure}[h]
\includegraphics{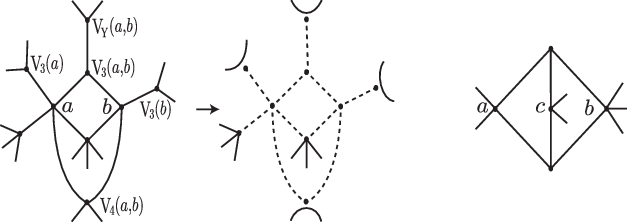}\\
\vspace{3mm}
\hspace{20mm}(a)\hspace{55mm} (b)
\caption{Deriving $\widehat{G}_{a,b}$} \label{fig:3}
\end{figure}

First consider the graph $G \setminus \{a,b\}$ for some distinct vertices $a$ and $b$.
In this graph each vertex of $V_3(a,b)$ has degree $1$.
And each vertex of $V_3(a), V_3(b)$ (not in $V_3(a,b)$) and $V_4(a,b)$ has degree $2$.
To derive $\widehat{G}_{a,b}$, we first delete all edges incident to $a$ and $b$ from $G$,
and then we also delete the remaining edges incident to $V_3(a,b)$,
and finally we contract one edge of the remaining pair of edges incident to each vertex of
$V_3(a)$, $V_3(b)$ (not in $V_3(a,b)$), $V_4(a,b)$ and $V_Y(a,b)$ as dotted lines in Figure~\ref{fig:3}$(a)$.
Thus we have the following equation counting the number of edges of $\widehat{G}_{a,b}$
which is called a {\em count equation\/};
\[ |\widehat{E}_{a,b}| = 21 - |E(a)\cup E(b)| - (|V_3(a)|+|V_3(b)|-|V_3(a,b)|+|V_4(a,b)|+|V_Y(a,b)|) \]

For short, $NE(a,b) = |E(a)\cup E(b)|$ and $NV_3(a,b) = |V_3(a)|+|V_3(b)|-|V_3(a,b)|$.
If $a$ and $b$ are adjacent vertices (i.e. dist$(a,b)=1$),
then all of $V_3(a,b), V_4(a,b)$ and $V_Y(a,b)$ are empty because $G$ is triangle-free.
Note that this manner of deriving $\widehat{G}_{a,b}$ must be handled in slightly different way
when there is a vertex $c$ in $V$ such that more than one vertex of $V(c)$ are contained
in $V_3(a,b)$ as in Figure~\ref{fig:3}$(b)$.
In this case we usually delete or contract more edges incident to $c$ even though $c$ is not in $V_Y(a,b)$.

A graph is $n$-{\em apex\/} if one can remove $n$ vertices from the graph to obtain a planar graph.
The following lemma gives an important condition for a graph not to be intrinsically knotted.

\begin{lemma}\label{lem:2-apex}
\cite{BBFFHL,OT}
If $G$ is $2-$apex, then $G$ is not intrinsically knotted.
\end{lemma}

The following two lemmas play an important role for $G$ to be $2$-apex.

\begin{lemma}\label{lem:planar8}
If $|\widehat{E}_{a,b}| \leq 8$, then $\widehat{G}_{a,b}$ is a planar graph.
Thus $G$ is not intrinsically knotted.
\end{lemma}

\begin{lemma}\label{lem:planar9}
If $|\widehat{E}_{a,b}|=9$, then $\widehat{G}_{a,b}$ is either a planar graph or homeomorphic to $K(3,3)$.
Furthermore if $\widehat{G}_{a,b}$ is not homeomorphic to $K(3,3)$,
then $G$ is not intrinsically knotted.
\end{lemma}

The graph $K(3,3)$ is a bipartite graph where each part has three vertices and 
each vertex is adjacent to every vertex in the opposite part, 
and so it is a triangle-free graph and every vertex has degree $3$.

To prove Theorem \ref{thm:main1}, we will show that any triangle-free graph with $21$ edges is eventually
either a $2$-apex or homeomorphic to one of $H_{12}$ or $C_{14}$.
Since intrinsically knotted graphs have at least $21$ edges \cite{JKM},
it is sufficient to consider simple and connected graphs having no vertex of degree $1$ or $2$.
Our process is constructing all possible such triangle-free graph $G$ with $21$ edges,
deleting two suitable vertices $a$ and $b$ of $G$, and then counting the number of edges of $\widehat{G}_{a,b}$.
If $\widehat{G}_{a,b}$ has $9$ edges or less,
we can use Lemma \ref{lem:planar8} or \ref{lem:planar9} in order to show that $G$ is not intrinsically knotted.
In the event that $\widehat{G}_{a,b}$ is not planar, we will show that $G$ is homeomorphic to $H_{12}$ or $C_{14}$.

Before describing the proof of Theorem \ref{thm:main1}, we introduce more notation.
Since $G$ is triangle-free, for any vertex $a$ of $G$, no two vertices in $V(a)$ are adjacent.
This means that $E(b)$ and $E(c)$ do not contain an edge in common
for any two distinct vertices $b$ and $c$ in $V(a)$.

\begin{itemize}
\item $E^2(a) = \underset{b \in V(a)}{\bigcup}E(b)$
\item $E \setminus E^2(a) = \{e_1(a), \cdots, e_{21-n}(a)\}$ \ \ \ if \ $|E^2(a)| = n < 21$
\end{itemize}
$e_i(a)$ is called  an {\em extra edge\/},
and two endpoints of the edge are denoted as $x_i(a)$ and $y_i(a)$ where $\mbox{deg}(x_i(a))\geq \mbox{deg}(y_i(a))$.

\begin{figure}[h]
\includegraphics{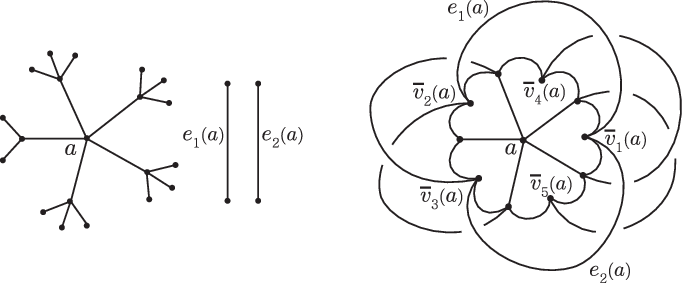}\\
\vspace{3mm}
(a)\hspace{60mm} (b)
\caption{Visualization of $G$}
\label{fig:4}
\end{figure}

In order to visualize $G$, we perform the following steps.
First choose a vertex $a$ with the maximal degree among all vertices and draw $E^2(a)$.
If $|E^2(a)|<21$, draw $E \setminus E^2(a)$
apart from $E^2(a)$ as in Figure~\ref{fig:4}$(a)$.
Then all vertices of degree $1$ of $E^2(a)$ and $E \setminus E^2(a)$
are merged into some vertices of degree at least $3$ without adding new edges as in Figure~\ref{fig:4}$(b)$.
Let $\overline{V}(a)$ denote the set of all such vertices, and let $[\overline{V}(a)]$ denote
a sequence of the degree of vertices in $\overline{V}(a)$ as follows:

\begin{itemize}
\item $\overline{V}(a) = V \setminus (V(a)\cup \{a\}) =
     \{\overline{v}_1(a), \cdots, \overline{v}_m(a)\}$ \
     with  $\mbox{deg}(\overline{v}_i(a))\geq \mbox{deg}(\overline{v}_{i+1}(a))$
\item $[\overline{V}(a)] =
     [\mbox{deg}(\overline{v}_1(a)), \cdots, \mbox{deg}(\overline{v}_m(a))]$
\item $|[\overline{V}(a)]| =
     \mbox{deg}(\overline{v}_1(a))+ \cdots + \mbox{deg}(\overline{v}_m(a))$
\end{itemize}
The graph in Figure~\ref{fig:4}$(b)$ is an example satisfying
$\mbox{deg}(a)=5,\ |V_3(a)|=1,\ |E^2(a)|=19$, and $[\overline{V}(a)]=[4,4,4,3,3]$.

The remaining three sections of the article are devoted to proof of Theorem \ref{thm:main1}.
From now on $a$ denotes one of vertices with maximal degree in $G$.
The proof is divided into three parts according to the degree of $a$.
In Section~\ref{sec:deg5} we show that any graph $G$ with $\mbox{deg}(a)\geq 5$ cannot be intrinsically knotted.
In Section~\ref{sec:deg4} we show that an intrinsically knotted graph with $\mbox{deg}(a)=4$ is exactly $H_{12}$.
Finally, in Section~\ref{sec:deg3} we show that any intrinsically knotted graph
all of whose vertices have degree $3$ is always $C_{14}$.

\section{$\deg(a) \geq 5$} \label{sec:deg5}

In this section we will show that for some $a', b' \in V$ either $|\widehat{E}_{a',b'}| \leq 8$
or $|\widehat{E}_{a',b'}|=9$, but that $\widehat{G}_{a',b'}$ is not homeomorphic to $K(3,3)$
by showing that it contains a vertex of degree more than $3$ or a triangle (or sometimes a bigon).
Then as a conclusion $G$ is not intrinsically knotted by Lemma~\ref{lem:planar8} and \ref{lem:planar9}.
Recall that $G$ has $21$ edges, every vertex has degree at least $3$, and $a$ has the maximal degree among them.

\subsection{Case $\deg(a) \geq 6$ or $\deg(a)=5$ with $|V_3(a)| \geq 4$}   \hspace{1cm}

If $\mbox{deg}(a) \geq 6$ then $|V_3(a)| \geq 3$.
Let $c$ be any vertex in $V_3(a)$.
Choose a vertex $b$ which has the maximal degree among $V(c) \setminus \{a\}$.
Then $|E(b)| + |V_Y(a,b)| \geq 4$ since $|V_Y(a,b)| \geq 1$ when $\mbox{deg}(b)=3$.
Note that $|V_3(b)| \geq |V_3(a,b)|$.
By the count equation, $|\widehat{E}_{a,b}| \leq 8$ in $\widehat{G}_{a,b}$.

Suppose that $\mbox{deg}(a)=5$ and $|V_3(a)|\geq 4$.
The proof is similar to the previous paragraph.

\subsection{Case $\deg(a)=5$ and $|V_3(a)|=3$}   \hspace{1cm}

Let $b$ and $c$ be two vertices of $V(a) \setminus V_3(a)$.
First, suppose that both of them have degree $5$.
Then $NE(a,b)=9$ and $|V_3(a)|=3$, so $|\widehat{E}_{a,b}| \leq 9$.
Furthermore, the vertex $c$ has degree $4$ in $\widehat{G}_{a,b}$, so it follows that
$\widehat{G}_{a,b}$ is not homeomorphic to $K(3,3)$.
Thus, $G$ is not intrinsically knotted by Lemma~\ref{lem:planar9}.

Now assume that one of them, say $b$, has degree $4$.
If $V(b) \setminus \{a\}$ consists of three vertices, all of which are of degree $3$,
then $NE(a,b)=8$ and $NV_3(a,b)=6$, so $|\widehat{E}_{a,b}| \leq 7$.
If not, let $d$ be a vertex of $V(b)$ which has degree at least $4$.
Then $NE(a,d) \geq 9$, $|V_3(a)|=3$, and $|V_4(a,d)| \geq 1$, because $V_4(a,d) \ni b$.
This implies that $|\widehat{E}_{a,d}| \leq 8$.

\subsection{Case $\deg(a)=5$ and $|V_3(a)|=0$} \hspace{1cm}

First, suppose that $V(a)$ contains a vertex of degree $5$, say $c$.
Since $G$ has $21$ edges, the other four vertices of $V(a)$ have degree $4$.
By the previous cases, it is sufficient to suppose that $|V_3(c)|\leq 2$.
So $V(c) \setminus \{a\}$ has at least two vertices, say $b$ and $d$, of degree $4$ or $5$.
Since $|E^2(a)|=21$ and $G$ is triangle-free,
all edges of $E(b)$ must be incident to different vertices of $V(a)$, so $|V_4(a,b)|\geq 3$.
This implies that $|\widehat{E}_{a,b}| \leq 9$.
Since $\widehat{G}_{a,b}$ has the vertex $d$ of degree at least $4$,
it follows that $\widehat{G}_{a,b}$ is not homeomorphic to $K(3,3)$.

Now, assume that all vertices of $V(a)$ have degree $4$ providing $|E^2(a)|=20$.
Let $e_1(a)$ be the extra edge and
recall that two endpoints of $e_1(a)$ are $x_1(a)$ and $y_1(a)$ with $\mbox{deg}(x_1(a))\geq \mbox{deg}(y_1(a))$.
Since $G$ is triangle-free,
all edges of $E(x_1(a)) \cup E(y_1(a))$ except $e_1(a)$ must be incident to different vertices of $V(a)$.
Thus the degrees of $x_1(a)$ and $y_1(a)$ must be either $4$ and $3$, or $3$ and $3$ respectively.
If $\mbox{deg}(x_1(a)) = 4$, then $|V_4(a,x_1(a))|=3$ and $|V_3(x_1(a))|=1$, so $|\widehat{E}_{a,x_1(a)}| = 8$.
If not, $[\overline{V}(a)]$ is either $[5,3,3,3,3]$ or $[4,4,3,3,3]$ because $|[\overline{V}(a)]|=17$.
Thus $\overline{v}_1(a)$ has degree $5$ or $4$ and differs from $x_1(a)$ and $y_1(a)$,
so $|V_4(a,\overline{v}_1(a))| \geq 4$.
Therefore, $|\widehat{E}_{a,\overline{v}_1(a)}| \leq 8$.

\subsection{Case $\deg(a)=5$ and $|V_3(a)|=1$} \hspace{1cm}

In this case, $V(a)$ contains four vertices of degree $4$ or $5$.
Let $n$ be the number of such vertices of degree $4$, and so we have $4-n$ of degree $5$, where $n=2,3,4$.
This implies that $|E^2(a)|=21 + (2-n)$ and $n-2$ extra edges exist.
If $\overline{V}(a)$ contains a vertex $\overline{v}_1(a)$ of degree $5$,
then five edges of $E(\overline{v}_1(a))$ are extra edges or incident to different vertices in $V(a)$.
For any of the above $n$, at least two among these edges are incident to vertices of degree $4$ in $V(a)$.
Then $NE(a,\overline{v}_1(a))=10$, $|V_3(a)|=1$, and $|V_4(a,\overline{v}_1(a))| \geq 2$,
implying $|\widehat{E}_{a,\overline{v}_1(a)}| \leq 8$.

Now, suppose that $\overline{V}(a)$ contains vertices of degree $3$ or $4$ only.
If $n=2$, $|[\overline{V}(a)]|=16$, and so $[\overline{V}(a)]$ is either $[4,4,4,4]$ or $[4,3,3,3,3]$.
For any vertex $b$ in $V_5(a)$,
four edges of $E(b)$ must be incident to different vertices of $\overline{V}(a)$.
Indeed these four edges are incident to four vertices of degree $4$, or
at least three edges among them are incident to vertices of degree $3$ in $\overline{V}(a)$.
This means that the vertex $b$ has degree $5$ with either $V_3(b)=0$ or $V_3(b) \geq 3$.
Both cases are dealt in previous cases 3.3, 3.1, and 3.2.

If $n=3$, $|[\overline{V}(a)]|=17$, and so $[\overline{V}(a)]=[4,4,3,3,3]$.
Let $V_5(a)=\{b\}$.
To avoid the case (3.2), four edges of $E(b)$ must be incident to
two vertices of degree $4$ and two vertices of degree $3$ in $\overline{V}(a)$,
which are $\overline{v}_1(a)$, $\overline{v}_2(a)$, $\overline{v}_3(a)$, and $\overline{v}_4(a)$.
Then there is a vertex $c$ of $V_4(a)$ such that at most one edge of $E(c)$ is
incident to $\overline{v}_3(a)$ and $\overline{v}_4(a)$,
ie two edges of $E(c)$ are incident to $\overline{v}_1(a)$, $\overline{v}_2(a)$, or $\overline{v}_5(a)$.
This implies that $NE(b,c)=9$ and $NV_3(b,c) + |V_4(b,c)| \geq 4$,
implying $|\widehat{E}_{b,c}| \leq 8$.

\begin{figure}[h]
\includegraphics{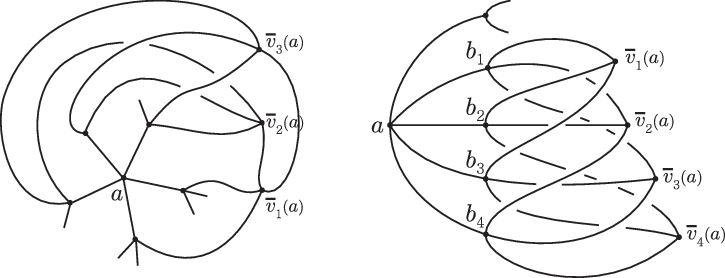}\\
\vspace{3mm}
(a)\hspace{60mm} (b)
\caption{$[4,4,4,3,3]$ and $[3,3,3,3,3,3]$ cases}
\label{fig:5}
\end{figure}

Finally, if $n=4$, $|[\overline{V}(a)]|=18$, and so $[\overline{V}(a)]$ is either $[4,4,4,3,3]$ or $[3,3,3,3,3,3]$.
Recall that two extra edges exist.
In the former case let $\{\overline{v}_1(a),\overline{v}_2(a),\overline{v}_3(a)\}$ be
the three vertices of degree $4$ in $\overline{V}(a)$.
For each $i=1,2,3$, if more than two edges of $E(\overline{v}_i(a))$ are incident to $V_4(a)$,
then $NE(a,\overline{v}_i(a))=9$, $|V_3(a)|=1$, and $|V_4(a,\overline{v}_i(a))| \geq 3$,
implying $|\widehat{E}_{a,\overline{v}_i(a)}| \leq 8$.
So, each of at least two edges of $E(\overline{v}_i(a))$ must be
either incident to the unique vertex of $V_3(a)$ or an extra edge.
Since $G$ is triangle-free, one of three vertices, say $\overline{v}_1(a)$,
has the property that $E(\overline{v}_1(a))$ contains both extra edges,
and $V(\overline{v}_1(a))$ and $V(\overline{v}_i(a))$ for each $i=2,3$
cannot share a vertex in $V(a)$.
This implies that $V(\overline{v}_2(a))$ and $V(\overline{v}_3(a))$ coincide
as in Figure~\ref{fig:5}$(a)$.
Then $NE(\overline{v}_2(a),\overline{v}_3(a))=8$, and
either $|V_4(\overline{v}_2(a),\overline{v}_3(a))|=4$ or
$|V_4(\overline{v}_2(a),\overline{v}_3(a))|=3$ and $|V_3(\overline{v}_2(a))|=1$.
Thus, $|\widehat{E}_{\overline{v}_2(a),\overline{v}_3(a)}| \leq 9$.
In $\widehat{G}_{\overline{v}_2(a),\overline{v}_3(a)}$ the vertex $a$ still has degree $4$ or $5$
so that $\widehat{G}_{\overline{v}_2(a),\overline{v}_3(a)}$ is not homeomorphic to $K(3,3)$.

In the latter case, let $V_4(a)=\{b_1,b_2,b_3,b_4\}$.
We claim that for some $i,j=1,2,3,4$, $|V_3(b_i,b_j)| \leq 1$.
Suppose not; that is, $|V_3(b_i,b_j)| \geq 2$ for all combinations of $i$ and $j$.
By some combinatorics we can derive that all $12$ edges of
$E(b_1) \cup E(b_2) \cup E(b_3) \cup E(b_4) \setminus E(a)$
are incident to only four vertices of $\overline{V}(a)$
as in Figure~\ref{fig:5}$(b)$.
This means that two extra edges must be incident to the remaining two vertices of $\overline{V}(a)$
at both endpoints.
But a bigon is not allowed.
Therefore without loss of generality $|V_3(b_1,b_2)| \leq 1$.
Then $NE(b_1,b_2)=8$ and $NV_3(b_1,b_2) \geq 5$,
implying $|\widehat{E}_{b_1,b_2}| \leq 8$.

\subsection{Case $\deg(a)=5$ and $|V_3(a)|=2$} \hspace{1cm}

If $V(a)$ contains a vertex of degree $5$, say $b$,
then previous four cases guarantee that we only consider that $|V_3(b)|=2$,
so $NV_3(a,b)=4$, which implies $|\widehat{E}_{a,b}|=8$.
Therefore we assume that $V(a)$ contains three vertices of degree $4$.
In this case three extra edges exist.
Since $|[\overline{V}(a)]|=19$, $[\overline{V}(a)]$ is
one of $[5,5,5,4]$, $[5,5,3,3,3]$, $[5,4,4,3,3]$, $[4,4,4,4,3]$, and $[4,3,3,3,3,3]$.

If, for some vertex $\overline{v}_i(a)$ with degree $5$,
one edge of $E(\overline{v}_i(a))$ is incident to $V_4(a)$,
then $NE(a,\overline{v}_i(a))=10$, $|V_3(a)|=2$, and $|V_4(a,\overline{v}_i(a))| \geq 1$,
implying $|\widehat{E}_{a,\overline{v}_i(a)}| \leq 8$.
Thus three edges of $E(\overline{v}_i(a))$ are extra edges and
the remaining two edges are incident to $V_3(a)$.
In the first two cases $[5,5,5,4]$ and $[5,5,3,3,3]$,
both $E(\overline{v}_1(a))$ and $E(\overline{v}_2(a))$ share three extra edges, but $G$ does not have a bigon.
In the third case $[5,4,4,3,3]$,
$E(\overline{v}_1(a))$ contains three extra edges and one of these extra edges must be incident to
$\overline{v}_4(a)$ or $\overline{v}_5(a)$ both of which have degree $3$.
Then $NE(a,\overline{v}_1(a))=10$ and $NV_3(a,\overline{v}_1(a)) \geq 3$,
implying $|\widehat{E}_{a,\overline{v}_1(a)}| \leq 8$.

If, for some vertex $\overline{v}_i(a)$ with degree $4$,
two edges of $E(\overline{v}_i(a))$ are incident to $V_4(a)$,
then $NE(a,\overline{v}_i(a))=9$, $|V_3(a)|=2$, and $|V_4(a,\overline{v}_i(a))| \geq 2$,
implying $|\widehat{E}_{a,\overline{v}_i(a)}| \leq 8$.
Thus at most one edge of $E(\overline{v}_i(a))$ is incident to $V_4(a)$.
In the fourth case, $[4,4,4,4,3]$,
at least twelve among sixteen edges incident to four vertices of degree $4$ in $\overline{V}(a)$
are not incident to $V_4(a)$.
This is impossible because there are only two vertices in $V_3(a)$ and three extra edges.
In the last case $[4,3,3,3,3,3]$,
since only one edge of $E(\overline{v}_1(a))$ is possibly incident to $V_4(a)$,
there is a vertex $b$ in $V_4(a)$ such that three edges of $E(b)$ are incident to
vertices of degree $3$ in $\overline{V}(a)$.
Then $NE(a,b)=8$ and $NV_3(a,b) \geq 5$, implying $|\widehat{E}_{a,b}| \leq 8$.

\section{$\deg(a) = 4$} \label{sec:deg4}

Since $|V|=|V_4|+|V_3|$ and $4|V_4| + 3|V_3| = 2|E|$,
the pair $(|V_4|,|V_3|)$ has three choices $(3,10)$, $(6,6)$ and $(9,2)$.
Here, $V_n$ denotes the set of vertices of degree $n$.
As in the preceding section, we will show that for some $a', b' \in V$ either $|\widehat{E}_{a',b'}| \leq 8$
or $|\widehat{E}_{a',b'}|=9$, but $\widehat{G}_{a',b'}$ is not homeomorphic to $K(3,3)$,
implying that $G$ is not intrinsically knotted.
But one exception occurs so that $G$ can possibly be $H_{12}$ when $(|V_4|,|V_3|)=(6,6)$.

\subsection{Case $(|V_4|,|V_3|)=(3,10)$} \hspace{1cm}

First suppose that $V_4$ has a vertex $a$ such that all four vertices of $V(a)$ have degree $3$.
Let $b_1$ and $b_2$ be the other vertices of $V_4$.
For each $i=1,2$, $NE(a,b_i)=8$.
If there is a vertex of $V_3(b_i)$ which is not contained in $V(a)$,
then $NV_3(a,b_i) \geq 5$, implying $|\widehat{E}_{a,b_i}| \leq 8$.
Thus each vertex of $V(b_1)$ is the vertex $b_2$ or contained in $V(a)$, and similarly for $b_2$.
This implies that the number of vertices of $V_3$ which have distance $1$ or $2$ from the vertex $a$
is at most $6$.
Take a vertex $c$ of $V_3$ with distance at least $3$ from $a$.
Since each vertex of $V(c)$ is neither $b_1$ nor $b_2$, it has degree $3$.
Thus $NE(a,c)=7$ and $NV_3(a,c) \geq 7$, implying $|\widehat{E}_{a,c}| \leq 7$.

Now, we only need to consider the case that 
each vertex of $V_4$ is adjacent to at least one vertex of degree $4$.
Then, without loss of generality, we have three vertices $a$, $b$, and $c$ of $V_4$
such that $V(b)$ contains $a$ and $c$.
If $V_3(a)$ and $V_3(c)$ do not coincide,
then $|V_4(a,c)|=1$ and $NV_3(a,c) \geq 4$, implying $|\widehat{E}_{a,c}| \leq 8$.
If $V_3(a)$ and $V_3(c)$ coincide and $|V_Y(a,c)| \geq 2$,
then $|V_4(a,c)|=1$ and $NV_3(a,c)=3$, implying $|\widehat{E}_{a,c}| \leq 7$.
If not, for the unique vertex $d$ of $V_Y(a,c)$, $V_3(a)=V_3(c)=V(d)$.
Then for a vertex $b'$ of $V_3(b)$, $V_3(b')$ is disjoint from $V_3(a)$.
Thus $NE(a,b')=7$, $NV_3(a,b')=5$, and $|V_4(a,b')|=1$, implying $|\widehat{E}_{a,b'}| \leq 8$.

\subsection{Case $(|V_4|,|V_3|)=(6,6)$} \hspace{1cm}

Consider the subgraph $H$ of $G$ consisting of all edges whose both end vertices have degree $4$.
Since $G$ has six vertices of degree $3$ and the same number of vertices of degree $4$,
$H$ is not empty set.
\vspace{2mm}

\noindent {\bf Claim 1.} \hspace{1mm}
If $H$ has a vertex of degree $1$, then $G$ is not intrinsically knotted.

\begin{figure}[h]
\includegraphics{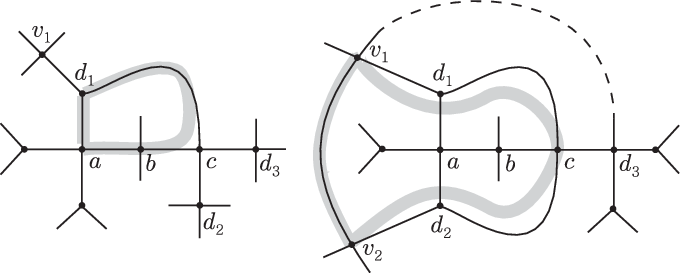}\\
\vspace{3mm}
\hspace{-11mm}  (a)\hspace{52mm} (b)
\caption{Some non-intrinsically knotted cases}
\label{fig:6}
\end{figure}

\begin{proof}
Suppose that $H$ has a vertex $a$ of degree $1$.
Let $b$ be the unique vertex of degree $4$ in $V(a)$.
If $|V_3(b)|=3$, then $NE(a,b)=7$ and $NV_3(a,b)=6$, implying $|\widehat{E}_{a,b}| \leq 8$.
Thus, there is another vertex $c$ of $V_4(b)$, and so we let $V(c)= \{b,d_1,d_2,d_3\}$.

First, assume that $|V_3(c)|=0$.
So the two vertices of $V(b) \setminus \{a,c\}$ must have degree $3$
because the six vertices $a, b, c, d_1, d_2$, and $d_3$ in $V_4$ are all different.
Thus $NE(a,b)=7$ and $NV_3(a,b)=5$, so $|\widehat{E}_{a,b}| \leq 9$.
Since $\widehat{G}_{a,b}$ has another vertex $d_1$ of degree $4$,
it follows that $\widehat{G}_{a,b}$ is not homeomorphic to $K(3,3)$.

Second, assume that $|V_3(c)|=1$, say $d_1 \in V_3(c)$.
If $d_1$ is not one of the vertices in $V(a)$, then $NE(a,c)=8$ and $NV_3(a,c) + |V_4(a,c)| = 5$, 
implying $|\widehat{E}_{a,c}| \leq 8$.
So we may assume that $d_1$ is in $V(a)$ and let $V(d_1) = \{a,c,v_1\}$.
If $v_1$ has degree $3$, then $NV_3(a,c) + |V_4(a,c)| = 4$ and $V_Y(a,c)=\{v_1\}$, implying $|\widehat{E}_{a,c}| \leq 8$.
Otherwise $v_1$ has degree $4$ and it is different from $d_2$ and $d_3$.
For any $i=2,3$, each vertex of $V(d_i) \setminus \{c\}$ either has degree $3$ or is $v_1$.
Thus $NE(d_2,d_3)=8$ and $NV_3(d_2,d_3) + |V_4(d_2,d_3)| \geq 4$, implying $|\widehat{E}_{d_2,d_3}| \leq 9$.
But $\widehat{G}_{d_2,d_3}$ has a triangle containing vertices $a$, $b$ and $d_1$.
See Figure~\ref{fig:6}(a).

Last, assume that $|V_3(c)| \geq 2$ and let $d_1$ and $d_2$ be two such vertices.
As in the previous case, we may say that $d_1$ and $d_2$ are in $V(a)$,
and $V(d_i) = \{a,c,v_i\}$ for $i=1,2$ where $v_i$ has degree $4$.
When $v_1=v_2$, $|V_3(a)|=3$, $|V_4(a,c)|=1$,
and $v_1$ has degree $2$ when we construct $\widehat{G}_{a,c}$,
implying $|\widehat{E}_{a,c}| \leq 8$.
When dist$(v_1,v_2) \geq 2$, three cases occur as follows:
$|V_3(v_1)| \geq 3$, $|V_3(v_2)| \geq 3$,
or for both $i=1,2$ $|V_3(v_i)| = 2$ and $V_4(v_i) = V_4 \setminus \{a,c,v_1,v_2\}$.
All three cases satisfy that $NV_3(v_1,v_2) + |V_4(v_1,v_2)| \geq 4$,
implying $|\widehat{E}_{v_1,v_2}| \leq 9$.
But $\widehat{G}_{v_1,v_2}$ has a bigon containing vertices $a$ and $c$.
Finally, when dist$(v_1,v_2)=1$, two cases occur as follows.
If $d_3$ has degree $3$, then by the same reason as before we may say that $d_3$ is also in $V(a)$,
and $V(d_3) = \{a,c,v_3\}$  where $v_3$ has degree $4$.
By the previous argument any pair of $v_1$, $v_2$ and $v_3$ has distance $1$.
This implies that $G$ contains a triangle.
If $d_3$ has degree $4$, then $|V_3(d_3)| \geq 2$ because at most one vertex of $V(d_3)$ can be $v_1$ or $v_2$.
Thus $NV_3(a,d_3) \geq 4$, implying $|\widehat{E}_{a,d_3}| \leq 9$.
But $\widehat{G}_{a,d_3}$ has a triangle containing vertices $c$, $v_1$ and $v_2$.
See Figure~\ref{fig:6}(b).
\end{proof}

\noindent {\bf Claim 2.} \hspace{1mm}
If $H$ is not a cycle with $6$ edges, then $G$ is not intrinsically knotted.

\begin{proof}
By Claim 1, if $H$ is not a cycle with $6$ edges,
then $H$ contains a cycle with $4$ or $5$ edges.
First assume that $H$ contains a cycle with $5$ edges.
Let $\{a_1,\cdots,a_5\}$ be the set of five vertices of the cycle appearing in clockwise order.
If the remaining vertex $b$ of $V_4$ is contained in some $V(a_i)$, say $i=1$,
then $b$ must have distance $1$ from one of $a_3$ and $a_4$, say $a_3$, by Claim 1.
See Figure~\ref{fig:7}.
If $V_3(a_2) \neq V_3(b)$, $NV_3(a_2,b) + |V_4(a_2,b)| \geq 5$, implying $|\widehat{E}_{a_2,b}| \leq 8$.
Otherwise, $V_3(a_2) = V_3(b)$.
Let $c_1$ and $c_3$ be the vertices of $V_3(a_1)$ and $V_3(a_3)$, respectively.
If $c_1=c_3$, we still have $|\widehat{E}_{a_2,b}| \leq 9$
and $\widehat{G}_{a_2,b}$ has a triangle containing vertices $a_5$, $a_4$ and $c_1=c_3$.
If $c_1 \neq c_3$, then $|\widehat{E}_{a_1,a_3}| \leq 9$
and $\widehat{G}_{a_1,a_3}$ has a bigon as in the figure.

If $b$ is not contained in $V(a_i)$ for any $i=1, \cdots, 5$, then $|V_3(a_i)|=2$.
If there is a pair of vertices $a_i$ and $a_{i+2}$ (or $a_{i-3}$ if $i=4,5$)
such that $V_3(a_i)$ and $V_3(a_{i+2})$ are disjoint,
then $NV_3(a_i,a_{i+2}) + |V_4(a_i,a_{i+2})| = 5$, implying $|\widehat{E}_{a_i,a_{i+2}}| \leq 8$.
Otherwise, for any pair of vertices $a_i$ and $a_{i+2}$ (or $a_{i-3}$ if $i=4,5$),
$V_3(a_i)$ and $V_3(a_{i+2})$ share vertices.
Then they must share only one vertex as in Figure~\ref{fig:7}(b).
Since there is only one extra vertex $b$ of degree $4$,
for some pair of vertices $a_i$ and $a_{i+2}$,
$NV_3(a_i,a_{i+2}) + |V_4(a_i,a_{i+2})| = 4$ and $V_Y(a_i,a_{i+2}) \geq 1$,
implying $|\widehat{E}_{a_i,a_{i+2}}| \leq 8$.

\begin{figure}[h]
\includegraphics{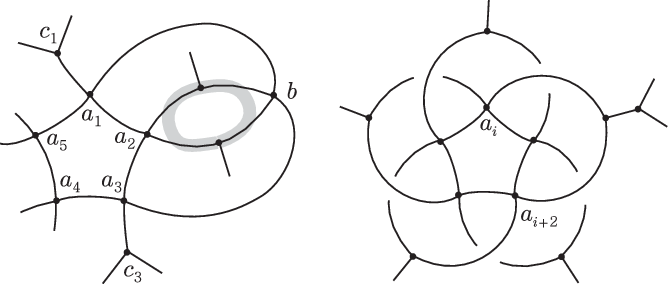}\\
\vspace{3mm}
\hspace{-7mm}  (a)\hspace{50mm} (b)
\caption{Cycle with $5$ edges}
\label{fig:7}
\end{figure}

Now, assume that $H$ contains a cycle with $4$ edges.
Let $\{a_1,\cdots,a_4\}$ be the set of four vertices of the cycle appearing in clockwise order.
If $V(a_1)$ and $V(a_3)$ (or similarly for $V(a_2)$ and $V(a_4)$) share only two vertices $a_2$ and $a_4$,
then the remaining two vertices of $V_4$ must be contained in $V(a_1) \cup V(a_3)$.
Otherwise, since $V(a_1) \cup V(a_3)$ has four more vertices other than $a_2$ and $a_4$, 
$NV_3(a_1,a_3) \geq 3$ and $|V_4(a_1,a_3)| = 2$, implying $|\widehat{E}_{a_1,a_3}| \leq 8$.
By Claim 1, the two vertices have distance $1$,
so $H$ contains a cycle with $5$ edges which was dealt in the previous case.
If $V(a_1)$ and $V(a_3)$ (or similarly for $V(a_2)$ and $V(a_4)$) share exactly three vertices, $a_2$, $a_4$ and $b$,
then let $c_1$ and $c_3$ be the remaining vertices of $V(a_1)$ and $V(a_3)$, respectively.
If both $c_1$ and $c_3$ have degree $3$, then $NV_3(a_1,a_3) + |V_4(a_1,a_3)| \geq 5$.
If both have degree $4$, then $H$ contains a cycle with $5$ edges as in the previous case.
Finally, if only $c_1$ (or similarly $c_3$) has degree $4$,
then, by Claim 1, $V(c_1)$ contains another vertex, say $d$, of $V_4$ and
also $d$ must have distance $1$ from one of $a_2$ and $a_4$, say $a_4$, as in Figure~\ref{fig:8}(a).
So $NV_3(a_4,c_1) + |V_4(a_4,c_1)| \geq 4$, implying $|\widehat{E}_{a_4,c_1}| \leq 9$,
and $\widehat{G}_{a_4,c_1}$ has a triangle containing vertices $a_2$, $a_3$, and $b$.
Now we may assume that $V(a_1) = V(a_3)$ and $V(a_2) = V(a_4)$.
Then $NV_3(a_1,a_3) + |V_4(a_1,a_3)| = 4$, implying $|\widehat{E}_{a_1,a_3}| \leq 9$,
and so $\widehat{G}_{a_1,a_3}$ has a bigon as in Figure~\ref{fig:8}(b).
\end{proof}

\begin{figure}[h]
\includegraphics{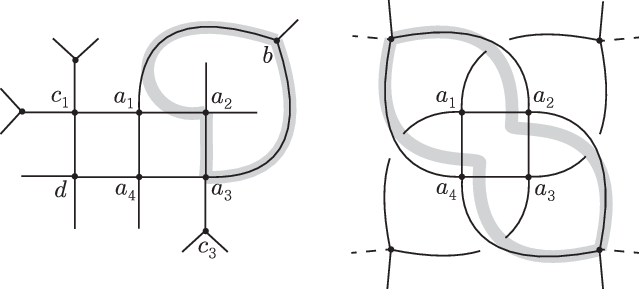}\\
\vspace{3mm}
\hspace{-3mm}  (a)\hspace{53mm} (b)
\caption{Cycle with $4$ edges}
\label{fig:8}
\end{figure}

\begin{figure}[h]
\includegraphics{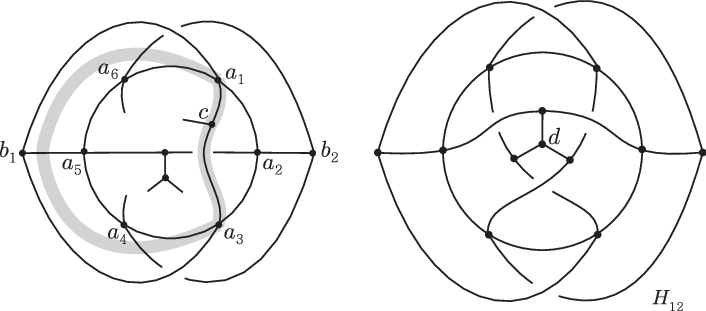}\\
\vspace{3mm}
(a)\hspace{58mm} (b) \\
\vspace{10mm}
\includegraphics{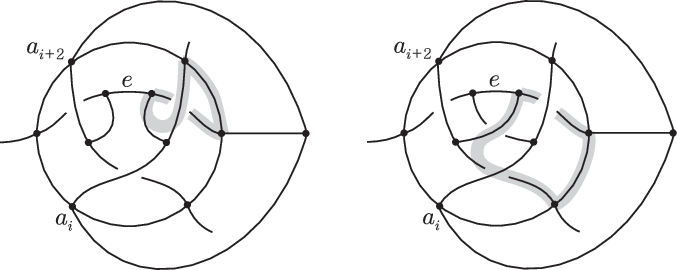}\\
\vspace{3mm}
(c)
\caption{Constructing $H_{12}$}
\label{fig:9}
\end{figure}

By Claim 2, $H$ is exactly a cycle with $6$ edges.
Let $\{a_1,\cdots,a_6\}$ be the set of six vertices of the cycle
with $a_i$ adjacent to $a_{i+1}$ for $i=1, \cdots, 5$, and $a_6$ adjacent to $a_1$.
First, suppose that there is not a vertex $b$ in $V_3$ such that $V(b)= \{a_1, a_3, a_5\}$.
If $V_3(a_1)$ and $V_3(a_3)$ are disjoint, then $NV_3(a_1,a_3) + |V_4(a_1,a_3)| = 5$.
If $V_3(a_1)$ and $V_3(a_3)$ share exactly one vertex $c$,
then the vertex of $V(c) \setminus \{a_1, a_3\}$ is not $a_5$, so it should be one of $V_Y(a_1,a_3)$.
Thus $NV_3(a_1,a_3) + |V_4(a_1,a_3)| + |V_Y(a_1,a_3)| = 5$.
If $V_3(a_1)$ and $V_3(a_3)$ are same, then $NV_3(a_1,a_5) + |V_4(a_1,a_5)| = 5$,
because $V_3(a_1)$ and $V_3(a_5)$ are disjoint.
All three cases guarantee that $G$ is not intrinsically knotted.
Therefore we may assume that there are two vertices $b_1$ and $b_2$
so that $V(b_1)= \{a_1, a_3, a_5\}$ and $V(b_2)= \{a_2, a_4, a_6\}$.
See Figure~\ref{fig:9}(a).

Suppose that there is a vertex $c$, with $c \neq b_1$, so that $V(c)$ contains $a_1$ and $a_3$.
Let $d_2$ and $d_5$ be the vertices of $V_3(a_2)$ and $V_3(a_5)$, other than $b_1$ and $b_2$, respectively.
If $d_2 \neq d_5$, then $NV_3(a_2,a_5) = 4$.
If $d_2 = d_5$, then $NV_3(a_2,a_5) = 3$ and $V_Y(a_2,a_5)$ is not empty.
Both cases provide $|\widehat{E}_{a_2,a_5}| \leq 9$,
and $\widehat{G}_{a_2,a_5}$ has a triangle containing vertices $a_1$, $a_3$ and $c$.
Therefore we may assume in general that for any vertex $c$, except $b_1$ and $b_2$,
$V(c)$ does not contain both $a_i$ and $a_{i+2}$ for any $i=1,2,3,4$,
and both $a_i$ and $a_{i-4}$ for any $i=5,6$.

Now, we conclude that $E \setminus \{ E^2(b_1) \cup E^2(b_2) \}$ consists of three extra edges.
Note that each vertex of these edges has degree $3$ and
there are four more vertices of degree $3$, besides $b_1$ and $b_2$.
These two facts guarantees that these extra edges must be connected as a tree.
This tree can be of two types; either all three edges are incident to one vertex $d$,
or two edges are incident to different endpoints of the other edge $e$, respectively.
In both cases, any two edges adjoined to the tree at the same vertex at the end
must be also incident to $a_i$ and $a_{i+3}$, respectively, for some $i=1,2,3$.
Therefore, $G$ is one of three graphs as in Figure~\ref{fig:9}(b)-(c), depending on the types of the tree.
The graph $G$ in Figure~\ref{fig:9}(b) is $H_{12}$, which is intrinsically knotted.
But the two graphs in Figure~\ref{fig:9}(c) are not intrinsically knotted
because, for some $i$, $|\widehat{E}_{a_i,a_{i+2}}| \leq 9$,
and $\widehat{G}_{a_i,a_{i+2}}$ has a triangle.

\subsection{Case $(|V_4|,|V_3|)=(9,2)$} \hspace{1cm}

Let $b_1$ and $b_2$ be the vertices of $V_3$.
Since $|V_3|=2$, there are at least three vertices, $a_1$, $a_2$, and $a_3$, in $V_4$
such that all vertices of each $V(a_i)$ have degree $4$.
If dist$(a_1,a_2)=1$, then $V(a_1) \cup V(a_2)$ consists of $8$ vertices of $V_4$,
and so let $c$ be the ninth vertex.
Let $d$ be any vertex among $V(a_1) \cup V(a_2) \setminus \{a_1,a_2\}$
which is not contained in $V(c)$.
We assume that $d$ is in $V(a_1)$.
Then $V(d)$ should be contained in $V(a_2) \cup \{b_1,b_2\}$.
This implies that $NE(a_2,d)=8$ and $|V_3(d)| + |V_4(a_2,d)| \geq 4$,
implying $|\widehat{E}_{a_2,d}| \leq 9$.
Since $c$ has degree $4$ in $\widehat{G}_{a_2,d}$,
it follows that $\widehat{G}_{a_2,d}$ is not homeomorphic to $K(3,3)$.
We have the same result for any choices of pairs among $a_1$, $a_2$, and $a_3$.

Now assume that the distance between any pair among $a_1$, $a_2$, and $a_3$ is at least $2$.
We separate into several cases according to the number $|V_4(a_1,a_2)|$.
If $V_4(a_1,a_2) = \emptyset$ (i.e. dist$(a_1,a_2) > 2$), then $|V_4| \geq 10$, a contradiction.
If $V_4(a_1,a_2)= \{d\}$, then $V_4 = V(a_1) \cup V(a_2) \cup \{a_1,a_2\}$.
This implies that $a_3 \in V(a_1) \cup V(a_2)$,
so dist$(a_1,a_3) = 1$ or dist$(a_2,a_3) = 1$
both of which were dealt with in the previous case.
If $V_4(a_1,a_2)= \{d_1,d_2\}$, then $V(d_1) \cup V(d_2) \setminus \{a_1,a_2\}$ is
contained in $\{a_3,b_1,b_2\}$.
This implies that each $V(d_i) \setminus \{a_1,a_2\}$ is a set of two vertices among $\{a_3,b_1,b_2\}$
so that $|V_3(d_1,d_2)| + |V_4(d_1,d_2)| \geq 4$, implying $|\widehat{E}_{d_1,d_2}| \leq 9$.
Since at least two of four vertices in $V(a_1) \cup V(a_2) \setminus \{d_1,d_2\}$
still have degree $4$ in $\widehat{G}_{d_1,d_2}$, it follows that $\widehat{G}_{d_1,d_2}$ is not homeomorphic to $K(3,3)$.
If $V_4(a_1,a_2) = \{d_1,d_2,d_3\}$,
then $V(d_1) \cup V(d_2) \cup V(d_3) \setminus \{a_1,a_2\}$ is contained in $\{a_3,a_4,b_1,b_2\}$,
where $a_3$ and $a_4$ are the remaining two vertices of degree $4$ other than $V(a_1) \cup V(a_2) \cup \{a_1,a_2\}$.
Thus each $V(d_i) \setminus \{a_1,a_2\}$ is the set of two vertices among $\{a_3,a_4,b_1,b_2\}$.
This implies that $|V_3(d_i,d_j)| + |V_4(d_i,d_j)| \geq 4$ for some $i,j=1,2,3$,
implying $|\widehat{E}_{d_i,d_j}| \leq 9$.
Since at least one of three vertices $V(a_1) \cup V(a_2) \setminus \{d_i,d_j\}$
still has degree $4$ in $\widehat{G}_{d_i,d_j}$, it follows that $\widehat{G}_{d_i,d_j}$ is not homeomorphic to $K(3,3)$.
Finally, if $|V_4(a_1,a_2)|=4$, then $|\widehat{E}_{a_1,a_2}| \leq 9$.
Since $\widehat{G}_{a_1,a_2}$ still has the remaining three vertices of degree $4$,
it follows that $\widehat{G}_{a_1,a_2}$ is not homeomorphic to $K(3,3)$.

\section{$\deg(a) = 3$} \label{sec:deg3}

\begin{figure}[h]
\includegraphics{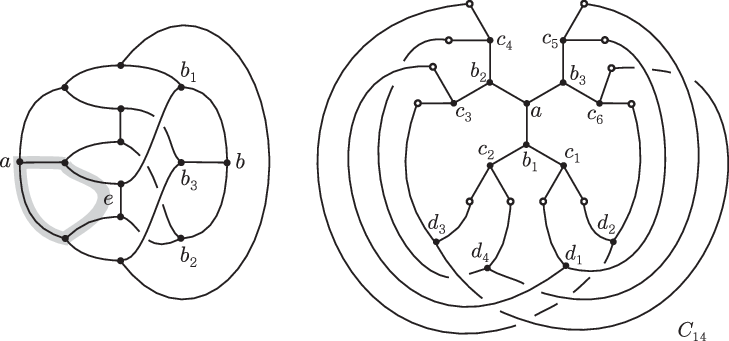}\\
\vspace{3mm}
(a)\hspace{58mm} (b)\caption{Constructing $C_{14}$}
\label{fig:10}
\end{figure}

Since we are working on the graph with $21$ edges and every vertex has degree $3$, there are exactly 14 vertices.
First, suppose that there exists a pair of vertices $a$ and $b$ with $\mbox{dist}(a,b) \geq 4$.
Then $E^2(a)$ and $E^2(b)$ can share vertices, but they do not share edges in common.
Since $|E^2(a) \cup E^2(b)| = 18$ and $|V(a) \cup V(b) \cup \{a,b\}| = 8$,
the $18$ endpoints of $E^2(a)$, $E^2(b)$ and three extra edges which are $E \setminus \{E^2(a) \cup E^2(b)\}$
meet at six vertices.
If any two edges of $E^2(a) \setminus E(a)$ (and similarly for $b$)
are incident to one vertex $c$ of these six vertices,
take the unique vertex $d$ of $V(a)$ which is not an endpoint of these two edges.
Then $NE(b,d)=6$ and $NV_3(b,d)=6$, implying $|\widehat{E}_{b,d}|=9$.
But $\widehat{G}_{b,d}$ has a triangle containing $c$ and the two vertices of $V(a) \setminus \{d\}$,
so it follows that $\widehat{G}_{b,d}$ is not homeomorphic to $K(3,3)$.
If not, each of these six vertices is a common endpoint of one edge of $E^2(a)$,
one edge of $E^2(b)$, and one extra edge.
Now take an extra edge $e$ and let $b_1$ and $b_2$ be the two vertices of $V(b)$
which have distance $1$ from the endpoints of $e$.
Let $b_3$ be the remaining vertex of $V(b)$.
Then $NE(b_1,b_2)=6$, $NV_3(b_1,b_2) = 5$, and $V_Y(b_1,b_2) = \{b_3\}$, implying $|\widehat{E}_{b_1,b_2}|=9$.
But $\widehat{G}_{b_1,b_2}$ has a triangle containing $a$ and two vertices of $V(a)$,
so it follows that $\widehat{G}_{b_1,b_2}$ is not homeomorphic to $K(3,3)$.
See Figure~\ref{fig:10}(a).

Therefore, we assume that the distance between any pair of vertices cannot exceed $3$.
Now we construct the intrinsically knotted graph $G$ satisfying these conditions.
Take a vertex $a$ and let $V(a)=\{b_1,b_2,b_3\}$ and $V(b_i)=\{a,c_{2i-1},c_{2i}\}$ for $i=1,2,3$.
As in Figure~\ref{fig:10}(b), the graph $E(a) \cup E(c_1) \cup \cdots \cup E(c_6)$
consists of $21$ edges and $22$ vertices.
We show this is the only way to draw the graph with $21$ edges
such that all vertices have distance at most 3 from $a$ and $10$ vertices
$a, b_1,b_2,b_3, c_1, \dots ,c_5$, and $c_6$ have degree $3$.
Now we join 12 white dots in Figure 10(b) into $4$ groups
indicating the remaining $4$ vertices $d_1,d_2,d_3$ and $d_4$.
Thus each $V(d_j)$, $j=1,2,3,4$, has three vertices among $c_1, \dots , c_6$.
Since the distance between any $c_i$ and $c_{i'}$ cannot exceed $3$,
the following two properties must be satisfied.
First property is that $V(d_j)$ contains exactly one vertex from each group $\{c_{2i-1},c_{2i}\}$ for $i=1,2,3$.
For example, if $V(d_1) = \{c_1,c_2,c_3\}$ (i.e. two vertices from the group $\{c_1,c_2\}$),
then we can connect $c_1$ to at most two vertices among $\{c_4,c_5,c_6\}$ through some $E(d_j)$.
This means that the distance between $c_1$ and one among $\{c_4,c_5,c_6\}$ exceed $3$.
The second property is that different $V(d_j)$ and $V(d_{j'})$ share at most one vertex.
For example, if they share two vertices $c_1$ and $c_3$, then $\mbox{dist}(c_1,c_4)=4$.
From these two properties, without loss of generality,
we may say that $V(d_1)=\{c_1,c_3,c_5\}$, $V(d_2)=\{c_1,c_4,c_6\}$,
$V(d_3)=\{c_2,c_3,c_6\}$, and $V(d_4)=\{c_2,c_4,c_5\}$ as drawn in Figure~\ref{fig:10}(b).
This graph is exactly $C_{14}$.

\end{document}